\documentclass[12pt]{amsart}
\usepackage{amsfonts,amssymb}
\usepackage[all]{xypic}
\usepackage{amsmath}
\usepackage{amsthm,amscd,latexsym}
\usepackage{color}

\newtheorem{theorem}{Theorem}[section]
\newtheorem{corollary}[theorem]{Corollary}
\newtheorem{Definition}[theorem]{Definition}

\newtheorem{proposition}[theorem]{Proposition}
\newtheorem{Example}[theorem]{Example}
\newtheorem{Remark}[theorem]{Remark}

\newenvironment{remark}{\begin{Remark}\begin{em}}{\end{em}\end{Remark}}

\newenvironment{definition}{\begin{Definition}\begin{em}}{\end{em}\end{Definition}}

\newcommand{\R}{{\mathbb R}}
\newcommand{\C}{{\mathbb C}}

\newcommand{\PP}{\mathbb P}
\newcommand{\N}{\mathbb N}

\newcommand{\Bx}{{\mathbf x}}
\newcommand{\By}{{\mathbf y}}

\newcommand{\Bw}{{\mathbf w}}

\DeclareMathOperator{\argmin}{\mathrm{arg\, min}}
\DeclareMathOperator{\Pro}{\mathcal{P}}

\begin{document}
\title[The Matrix Geometric Mean]{The Expanding Universe of the Geometric Mean}
\author{Jimmie D. Lawson and Yongdo Lim}

\address{Department of Mathematics, Louisiana State University,
Baton Rouge, LA70803, USA}\email{lawson@math.lsu.edu}

\address{ Department of Mathematics,
Sungkyunkwan University, Suwon 440-746, Korea} \email{ylim@skku.edu}

\keywords{geometric mean, operator mean,  positive definite matrices, Hadamard spaces,
contractive barycentric map}
\date{\today}




\keywords{geometric mean, operator mean,  positive definite matrices, Hadamard spaces,
contractive barycentric map}
\date{\today}
\maketitle
\begin{abstract}
In this paper the authors seek to trace in an accessible fashion the rapid recent development of the
theory of the matrix geometric mean in the cone of positive definite matrices up through
the closely related operator geometric mean in the positive cone of a unital $C^*$-algebra. The story begins with
the two-variable matrix geometric mean, moves to the breakthrough developments in the
multivariable matrix setting, the main focus of the paper, and then on to the extension to the positive
cone of the $C^*$-algebra of operators on a Hilbert space, even to general unital $C^*$-algebras, and
finally to the consideration of barycentric maps that grow out of the geometric mean  on the space
of integrable probability measures on the positive cone. Besides expected tools from linear algebra and operator theory,
one observes a surprisingly substantial interplay with geometrical notions in metric spaces, particularly the notion of nonpositive curvature.
Added features include a glance at the probabilistic theory of random variables with values in a metric space of nonpositive curvature,
and the appearance of related means such as the inductive and power means.
\end{abstract}

\section{Forward}
This manuscript is an expanded and  more detailed ArXiv version of the authors' earlier rather similar article \emph{Following the Trail
of the Operator Geometric Mean} \cite{LL20}.    Some twenty years ago the authors published an article \cite{LL01} in the \textsc{Monthly}
 that treated in some depth the two-variable matrix geometric mean and that has been rather widely read and cited.  In this paper
 we seek to describe the significant advancement and broad generalization of the theory that has taken place in these twenty years.

\section{Introduction}
The problem of ``squaring a rectangle" is the problem of constructing the side of a square that has the same area as a given rectangle.
Such a construction is given by Euclid in Book II of the \emph{Elements}.  If the sides of the rectangle are $a$ and $b$, then the side of the
square has length $\sqrt{ab}$, the geometric mean of $a$ and $b$.  Because of their interest in proportions and musical ratios,
the Greeks defined some 2500 years ago at least eleven different means, with the best known ones being the arithmetic, geometric, harmonic, and golden.
The subject of (binary) means for positive numbers or line segments has a rich mathematical lineage dating back into antiquity.
 The study of various means on the positive reals and their properties has continued off and on throughout the history of mathematics
 up to the present day.

The appropriate  definition of the geometric mean for two positive definite matrices of the same
size seems to have first appeared in 1975 in a paper of Pusz and Woronowicz \cite{PW}.
Ando \cite{An79} provided the first systematic development of many of
its basic properties and gave equivalent characterizations and also applications to matrix
inequalities that are otherwise difficult to prove.

In an article \cite{LL01} appearing some twenty years ago in the \textsc{Monthly}
the authors presented eight characterizations or prominent properties of the
classical geometric mean and showed how each extended to the matrix geometric mean setting,
providing convincing documentation that the name ``matrix geometric mean"  was
most appropriate.  As hinted in the Forward, that theory has now advanced to a
multivariable setting for both positive matrices and operators and beyond, as we will
trace out.

Positive definite matrices  have become fundamental computational
objects in many areas of engineering, statistics, quantum
information, applied mathematics, and elsewhere. They appear as ``data points" in
a diverse variety of settings: covariance matrices in statistics, elements of the search space in convex
and semidefinite programming, kernels in machine learning, observations in radar imaging,
and diffusion tensors in medical imaging, to cite only a few. A variety of computational algorithms
have arisen for approximation, interpolation, filtering, estimation, and averaging. Our interest
focuses on the last named, the process of finding an average or mean, which is again positive
definite. In recent years it has been increasingly recognized that the Euclidean distance is often not the most
suitable for the space $\PP$ of positive definite matrices and that working with the
appropriate geometry does matter in computational problems; see e.g. \cite{Bar}, \cite{MZ}.
The matrix geometric mean grows out of the geometric structure of $\PP$, which makes it
a particularly suitable averaging tool in a variety of settings.

\section{Positive Definite Matrices}
Let $\mathcal{M}_m(\C)$, or simply $\mathcal{M}_m$, denote the set of $m\times m$ complex matrices.  We may identify
$\mathcal{M}_m$ with the set of linear operators on $\C^m$, where we consider $\C^m$ to be a complex Hilbert space
of column vectors with the usual hermitian inner product. Denoting the conjugate transpose of $A\in\mathcal{M}_m$ by
$A^*$, we recall that $A$ is \emph{hermitian} if $A=A^*$ and \emph{unitary} if $A^*=A^{-1}$.  The hermitian matrix $A$ is
\emph{positive definite} if $\forall u\ne 0,\,\langle u,Au\rangle>0$.
These notions readily generalize to $\mathcal{B}(H)$, the algebra of operators on
an arbitrary Hilbert space.

The following are well-known equivalences for a hermitian matrix $A$
to be positive definite (with the definition appearing first):
\begin{enumerate}
\item $\langle Ax,x\rangle >0$ for all $0\ne x,$ where $\langle\cdot,\cdot\rangle$
is the Hilbert space inner product on ${\mathbb{C}}^m$.
\item  $A=BB^*$  for some invertible $B$.
\item  $A$ has all positive eigenvalues.
\item $A=\exp B=\sum_{k=0}^\infty B^k/k!$ for some (unique) hermitian $B$.
\item $A=UDU^*$ for some unitary $U$ and diagonal $D$ with positive diagonal entries.
\end{enumerate}

The positive definite  $m\times m$-matrices form an open cone in $\mathbb{H}_m$, the $m\times m$ hermitian
matrices, with closure
the positive semidefinite matrices (equivalently, $\langle Ax,x\rangle\geq 0$ for all
$x$).  We denote the open cone of positive definite matrices by $\PP$
(or $\PP_m$ if we need to distinguish the dimension). The exponential map
$\exp:\mathbb{H}\to \PP$ is an analytic diffeomorphism with inverse analytic
diffeomorphism $\log: \PP\to \mathbb{H}$.

Every positive definite (resp.\,hermitian) matrix operator has a unique \emph{spectral decomposition}
$$A=\sum_{i=1}^k\lambda_iE_i,$$
where the $\lambda_i>0$ (resp.\,$\lambda_i\in\R)$ range over the
distinct eigenvalues of $A$ and $E_i$ is the orthogonal projection
onto the eigenspace of $\lambda_i$.  One then has
$$A^p=\sum_{i=1}^k \lambda_i^p E_i,$$
from which one can easily deduce that every positive definite matrix
has a unique positive definite $p^{th}$-root. We also note that the
exponential map is given alternatively by
$$\exp A=\exp\left(\sum_{i=1}^k \lambda_i E_i\right)=\sum_{i=1}^k e^{\lambda_i}E_i,$$
from which we can quickly deduce the equivalence of items (3) and (4) in the previous list of equivalent characterizations
of positive definite matrices.

We define a partial order (sometimes called the \emph{Loewner order})
on the vector space $\mathbb{H}_m$ of hermitian matrices by $A\leq B$
if $B-A$ is positive semidefinite.  We note $0\leq A$ iff $A$ is
positive semidefinite and write $0<A$ if $A\in\PP$ iff $A$ is positive definite.
The matrix $A$ is sometimes called \emph{strictly positive} in this setting.

For any invertible $M\in \mathcal{M}_m(\mathbb{C})$, the
\emph{congruence transformation} $\Gamma_M(X)=MXM^*$  is an
invertible linear map on $\mathcal{M}_m$ that carries each of
$\mathbb{H}$, $\PP$, and $\overline{\PP}$, the convex cone of
positive semidefinite matrices, onto itself. It follows that
congruence transformations preserve the Loewner order on
$\mathbb{H}$:  $A\leq B$ implies $MAM^* \leq MBM^*$ for $M$
invertible.  We note, in particular, for each $M\in\PP$, $\Gamma_M(X)=MXM$ is
a congruence transformation on $\PP$.   Matrix
inversion $A\mapsto A^{-1}$ maps ${\Bbb P}$ onto itself and reverses
the  Loewner order, as we shall see later.

The geometry of $\PP$ will be crucial in what follows.  One
important approach to geometry is that of Felix Klein's Erlangen
Program, which emphasized the importance of the group of
transformations or ``symmetries"  that preserved basic geometric
properties. For the study of $\PP$ this group, denoted $G(\PP)$, is
the one generated by the congruence transformations and the
inversion map, which acts as the point reflection through
the identity matrix.
\begin{remark}\label{R:useful}
For $A,B\in\PP$, by the previous item (5), there exists a unitary $U$ such that $U(A^{-1/2}BA^{-1/2})U^*=D$, a diagonal matrix,
and hence $\Gamma_{UA^{-1/2}}$ carries $A$ to $I$ and $B$ to $D$.  This observation allows various results
about $A,B\in \PP$ to be reduced to the case $A=I$ and $B$ is a diagonal matrix.
\end{remark}

The arithmetic and harmonic means readily extend from
$\mathbb{R}^{>0}:=(0,\infty)$ to the set $\PP$ of
positive definite matrices:
$$\mathcal{A}(A,B)=\frac{1}{2}(A+B);  ~~~~~~\mathcal{H}(A,B)=2(A^{-1}+B^{-1})^{-1}.$$
The geometric mean is not so obvious (e.g., $\sqrt{AB},$
the square root of $AB$ with positive
eigenvalues, need not be positive definite for $A,B$ positive
definite). One approach is to rewrite the equation $x^2=ab$ (which
has positive solution the geometric mean of $a$ and $b$) in its
appropriate form in the noncommutative setting and solve for $X$:
\begin{eqnarray*}
XA^{-1}X&=&B\\
A^{-1/2}XA^{-1/2}A^{-1/2}XA^{-1/2}&=&A^{-1/2}BA^{-1/2}\\
A^{-1/2}XA^{-1/2} &=&(A^{-1/2}BA^{-1/2})^{1/2}\\
X &=& A^{1/2}(A^{-1/2}BA^{-1/2})^{1/2}A^{1/2}.
\end{eqnarray*}
\begin{definition}
The matrix geometric mean $A\# B$ of $A,B\in\PP$ is given by
$$A\# B=A^{1/2}(A^{-1/2}BA^{-1/2})^{1/2}A^{1/2}.$$
Alternatively it can be characterized as the unique positive definite solution $X$ of
the elementary Riccati equation $XA^{-1}X=B$.
\end{definition}

By inverting both sides $XA^{-1}X=B$ and multiplying through by $X$
on the right and left, one obtains $XB^{-1}X=A$.  Since the second
equation is equivalent to the first, we see $A\#B=B\#A$.  Similarly
one can use the Riccati equation to show that the matrix geometric
mean is invariant under congruence mappings and inversion, i.e.,
$M(A\#B)M^* =(MAM^*)\# (MBM^*)$ for $M$ invertible and
$(A\#B)^{-1}=A^{-1}\#B^{-1}.$ By the
Riccati equation $(A\#B)(A^{-1}-B^{-1})(A\#B)=B-A$, and this along with
order invariance under congruence maps shows that
the inversion map is order reversing. We collect these
properties together with other basic properties that can be deduced
by elementary arguments; see \cite{LL01}.
\begin{proposition}\label{P:basicprop}
The following hold in $\PP$:
\begin{itemize}
\item[(i)] $($Commutativity$)$ $A\#B=B\# A$.
\item[(ii)] $($Congruence Invariance$)$ $M(A\# B)M^*=MAM^*\#MBM^*$ for $M$ invertible.
\item[(iii)] $($Inversion Invariance$)$ $(A\# B)^{-1}=A^{-1}\#B^{-1}$.
\item[(iv)] $($Monotonicity$)$ If $C\leq A$ and $D\leq B$, then $C\#D\leq A\#B$.
\item[(v)]  For $AB=BA$, $A\#B=A^{1/2}B^{1/2}=(AB)^{1/2}$; in particular $A\#I=A^{1/2}$.
\item[(vi)] $($AGM Inequality$)$ $2(A^{-1}+B^{-1})^{-1}\leq A\#B\leq (1/2)(A+B)$.
\end{itemize}
\end{proposition}
 By the Loewner-Heinz inequality $A\#I=A^{1/2}\leq B^{1/2}=B\#I$ for $0<A\leq B$.
From this, using congruence invariance, one obtains the important property (iv).
As  mentioned in the introduction, other formulations and connections between
the matrix geometric mean and the one for positive real numbers may be found in \cite{LL01}.

\section{The Riemannian Metric}
The tools needed for extending the binary geometric mean on $\PP$ to a multivariable one have relied on the
metric and geometric structure of $\PP$.  Such considerations had already begun in the binary setting;
see \cite[Section 4]{LL01}.  We briefly overview the necessary tools in this section; see \cite[Section 4]{LL01} and
\cite[Chapter 6]{Bh} for details.

We  equip the space $\mathbb{H}$ of hermitian matrices of some fixed dimension $m$ with
the \emph{Frobenius} inner product $\langle A,B\rangle=$
Tr$(AB)$, the trace of $AB$, which makes $\mathbb{H}$ a Hilbert space. The corresponding norm $\Vert A\Vert_2=\langle A,A\rangle ^{1/2}$ is called the
\emph{Frobenius} or \emph{Hilbert-Schmidt} norm.  We can write $A=UDU^*$ for some unitary $U$, where $D$ is a diagonal matrix with entries
the eigenvalues of $A$.  We now observe that
\begin{equation}\label{E:Frob}
\Vert
A\Vert_2=(\mathrm{Tr}(UDU^*UDU^*))^{1/2}=(\mathrm{Tr}D^2)^{1/2}=\left(\sum_{i=1}^m
\lambda_i^2(A)\right)^{1/2},
\end{equation}
where $\{\lambda_i\}_{i=1}^m$ is the set of
eigenvalues of $A$.

We define the \emph{Riemannian metric} $\delta$ on $\PP$ by
\begin{equation}\label{E:Rmetric}
\delta(A,B)=\Vert \log(A^{-1/2}BA^{-1/2})\Vert_2=\left(\sum_{i=1}^m
\log^2\lambda_i(A^{-1/2}BA^{-1/2})\right)^{1/2},
\end{equation}
where $\{\lambda_i\}$ are the eigenvalues of $A^{-1/2}BA^{-1/2}$.  In the last expression we may replace $A^{-1/2}BA^{-1/2}$
by $BA^{-1}$, since the two are similar and hence have the same eigenvalues.

We list basic properties of the Riemannian metric $\delta$.
\begin{proposition}\label{P:metprops}
The Riemannian metric $\delta$ is a metric making $\PP$ a complete
metric space exhibiting  the following properties:
\begin{itemize}
\item[(1)] For $M\in \mathbf{GL}_m(\mathbb{C})$, the group of $m\times m$ invertible matrices, the
congruence transformation $\Gamma_M:(\PP,\delta)\to(\PP,\delta)$ defined by $\Gamma_M(X)=MXM^*$ is
an isometry.  The inversion map $A\mapsto A^{-1}$ is also an isometry on $\PP$.
\item[(2)] The exponential map $\exp: \mathbb{H}\to \PP$ is expansive, that is, $\delta(e^A,e^B)\geq \Vert B-A\Vert_2$
for all $A,B\in\mathbb{H}$.
\item[(3)] The exponential map restricted to any one-dimensional subspace of $\mathbb{H}$ is an isometry.  Furthermore,
any metric on $\PP$ that has this property and is invariant under congruence transformations must agree with $\delta$.
\item[(4)] For $A,B\in\PP$, $A\#B$ is the unique metric midpoint between $A$ and $B$.
\end{itemize}
\end{proposition}

\begin{proof} To give some flavor of the preceding, we prove parts of (1), (3) and (4).
First for (1) we observe for $A,B\in \PP$ and $M$ invertible, $(MBM^*)(MAM^*)^{-1}= M(BA^{-1})M^{-1}$, which
is similar to $BA^{-1}$, and hence to $A^{-1/2}BA^{-1/2}$, and thus may replace it in last part of equation
(\ref{E:Rmetric}).

For (3), let $A\in \mathbb{H}$.  Then
\begin{eqnarray*}
\delta(e^{sA}, e^{tA})&=&\delta(e^{(-s/2)A}e^{sA}e^{(-s/2)A},e^{(-s/2)A}e^{tA}e^{(-s/2)A})\\
&=&\delta(I,e^{(t-s)A})
=\Vert \log(e^{(t-s)A})\Vert_2=\Vert tA-sA\Vert_2. \end{eqnarray*}
Now let $d(\cdot,\cdot)$ be a metric satisfying the two properties of (3).
For $A,B\in\PP$,
\begin{eqnarray*}
 d(A,B)&=&d(A^{-1/2}AA^{-1/2}, A^{-1/2}BA^{-1/2})\\
 &=&d(I,A^{-1/2}BA^{-1/2})\\
 &=&\Vert \log(A^{-1/2}BA^{-1/2})\Vert _2=\delta(A,B).\end{eqnarray*}

For the midpoint property in (4), since congruence mappings preserve
$\delta$ and $\#$, it suffices to consider the case $A\#B=I$
(otherwise first apply $\Gamma_{(A\#B)^{-1/2}}$). Then
$B=IA^{-1}I=A^{-1}$.  Applying $\log$ we obtain $\log B=-\log A$, so
restricting $\exp$ to the one-dimensional subspace $\R\cdot A$ yields the result by (3).

A proof of (2) is sketched in \cite{LL01}, and a shorter and more
elegant proof appears in \cite[Chapter 6]{Bh}.

For the triangular inequality, consider first the case $A,B,C$ with $A=I$.
Then from (3) $\delta(I,B)=\Vert \log B\Vert_2$, $\delta(I,C)=\Vert \log C\Vert_2$ and
from (2) $\Vert\log(C)-\log(B)\Vert_2\leq \delta(B,C)$, so
$$\delta(I,C)=\Vert \log C\Vert_2\leq \Vert\log B\Vert_2+\Vert\log(C)-\log(B)\Vert_2\leq \delta(I,B)+\delta(B,C).$$
The general case follows using (1) and the congruence transformation $X\mapsto A^{-1/2}XA^{-1/2}$ and its inverse.

By (2) and continuity of the exponential map, the
metric $\delta$ is complete.
\end{proof}

For a metric space $(X,d)$ the metric is said to satisfy the \emph{semiparallelogram law} if
for all $x_1,x_2\in X$, there exists $m\in X$ such that for any $x\in X$,
$$d^2(x_1,x_2)+4d^2(x,m)\leq 2d^2(x,x_1)+2d^2(x,x_2)\qquad \qquad (NPC) $$
One can show that $m=m(x_1,x_2)$ is unique and is the unique metric midpoint between
$x_1$ and $x_2$. If the inequality is replaced by an equality, one obtains a version of the
parallelogram law holding in Hilbert spaces.  The semiparallelogram law is a metric version
of nonpositive curvature (NPC). We define a \emph{Hadamard space} to be a complete metric
space satisfying the semiparallelogram law. These spaces have been and continue to be
widely studied and appear under the alternative names of global NPC-spaces or  CAT(0)-spaces.

Using Proposition \ref{P:metprops} it is straightforward to show that $(\PP,\delta)$ is a
Hadamard space.  One first considers the case that $A\#B=I$ and another point $C$. Then the
parallelogram law holds in the Hilbert space $\mathbb{H}$ for $\log A$, $\log B=-\log A$, and $\log C$,
and one uses Proposition \ref{P:metprops}(2) to obtain the semiparallelogram law for $A,B,C$.
Via \ref{P:metprops}(1) the general case reduces to this one. See, for example,
\cite{LL01} for further details.
\begin{corollary}\label{C:NPC}
The space $(\PP,\delta)$ is a Hadamard space.
\end{corollary}

\begin{remark}\label{R:Riemann}
A more sophisticated approach to the results of this section is the path of Riemannian geometry.
The open cone $\PP_m$ of $m\times m$ positive
definite matrices becomes a well-known  Riemannian manifold when equipped with
the trace Riemannian metric: $\langle X,Y\rangle_A
={\mathrm{tr}} A^{-1}XA^{-1}Y,$ where $A\in \PP_m$
and $X,Y$ are $m\times m$ Hermitian matrices.  The corresponding
distance metric on $\PP_m$ is precisely our metric $\delta$, and
this is the source of the name ``Riemannian metric."
 The distance metric  of a simply connected Riemannian manifold satisfies (NPC) iff the manifold
has nonpositive curvature in the usual sense.
\end{remark}

\section{Means of several variables}
Formally a \emph{mean of order} $n$, or $n$-\emph{mean} for short, on a set $X$ is a function
$\gamma:X^n \to X$ satisfying the idempotency condition:  $\forall x\in X,\, \gamma(x,x,\ldots, x)=x.$
It is frequently assumed in the definition of a
mean that it is \emph{symmetric}, that is, invariant under any permutation of variables.
When we speak of an \emph{omnivariable} mean $\gamma=\{\gamma_n\},$ we are
referring to one defined for all $X^n$, $n \geq 1$.   (For $n=1$, $\gamma_1$ is the identity, and
is thus frequently ignored.) The mean $\gamma:X^n\to X$ is \emph{continuous}
or a \emph{topological mean} if $X$ is a topological space and $\gamma$ is continuous.  Frequently a mean represents some type of averaging operator.

In 2004 T. Ando, C. K. Li and R. Mathias \cite{ALM04} gave the first
extension of the binary geometric mean to an omnivariable mean on
$\PP$, which came to be called the ALM mean. For three variables
$A,B,C$, they first took the new three point set
$\{A_1:=B\#C, B_1:=A\#C, C_1:=A\#B\}$ consisting
of the geometric mean (midpoint) of each pair, then repeated this
construction on the new three point set and continued repeating the
operation inductively.  They showed the triples approached a common
point, their mean for the case $n=3$. They extended it inductively
to an $n$-variable mean for all $n>2$. We note that Lawson and Lim
\cite{LL08} later extended the ALM construction to a rather  wide
class of metric spaces.

Ando, Li, and Mathias made two important contributions in their
paper.  First of all, they identified axiomatic properties that an
omnivariable geometric mean $\gamma$ should satisfy. They then
established that the ALM mean they had defined satisfied all these
properties. The proofs typically involved extending from the known
case of $n=2$ by induction.

\begin{theorem}\label{T:ALM}
Let $\mathbb{A}=(A_1,\ldots, A_n),\mathbb{B}=(B_1,\ldots,B_n)\in\PP^n$ and
let $\gamma$ be the ALM mean. The following properties hold.
\begin{itemize}
\item[(P1)]$($Consistency with scalars$)$
$ \gamma({\mathbb A})=(A_{1}\cdots A_{n})^{1/n}$ if the $A_{i}$'s commute;
\item[(P2)] $($Joint homogeneity$)$
$ \gamma(a_{1}A_{1},\dots,a_{n}A_{n})= (a_{1}\cdots
a_{n})^{1/n}\gamma({\Bbb A});$

\item[(P3)] $($Permutation invariance$)$
$\gamma({\Bbb A}_{\sigma}) =\gamma({\Bbb A}),$ where ${\Bbb
A}_{\sigma}=(A_{\sigma(1)},\dots,A_{\sigma(n)});$

\item[(P4)] $($Monotonicity$)$ If $A_{i}\leq B_{i}$ for all $1\leq i\leq n,$ then
$ \gamma({\Bbb A})\leq \gamma({\Bbb B});$
\item[(P5)] $($Continuity$)$ $\gamma$ is continuous;

\item[(P6)] $($Congruence invariance$)$
$\gamma(M{\Bbb A}M^*)= M\gamma({\Bbb A})M^*$ for $M$ invertible,
where $M(A_{1},\dots,A_{n})M^{*}=(MA_{1}M^{*},\dots,MA_{n}M^{*});$

\item[(P7)] $($Joint concavity$)$
$\gamma(\lambda {\Bbb A}+(1-\lambda){\Bbb B})\geq \lambda \gamma({\Bbb
A})+(1-\lambda)\gamma({\Bbb B})$ for $0\leq \lambda \leq 1$;

\item[(P8)] $($Self-duality$)$
$\gamma(A_{1}^{-1} ,\dots, A_{n}^{-1})^{-1}= \gamma(A_{1} ,\dots, A_{n});$

\item[(P9)]$($Determinantal identity$)$
$ {\mathrm {Det}}\,\gamma({\Bbb A})= \prod_{i=1}^{n}({\mathrm
{Det}}A_{i})^{1/n}$; and

\item[(P10)] $($AGH mean inequalities$)$
$ n(\sum_{i=1}^{n}A_{i}^{-1})^{-1}\leq \gamma({\Bbb A}) \leq
\frac{1}{n}\sum_{i=1}^{n}A_{i}.$\\
$($AGH is short for arithmetic-geometric-harmonic.$)$
\end{itemize}
\end{theorem}
D. Bini, B. Meini, and F. Poloni \cite{BMP} later gave a variant of the ALM mean that retained its
properties, but was much more computationally efficient.

A  \emph{weight} $\Bw$ of length $n$ is an $n$-tuple $(w_1,\ldots,w_n)$ where $0< w_i\leq 1$ for each $i$
and $\sum_{k=1}^n w_k=1$, and a \emph{weighted $n$-tuple} of a set $X$ is a pair $(\Bw,\Bx)$, where
$\Bw$ is a weight of length $n$ and $\Bx=(x_1,\ldots,x_n)\in X^n$.  We think of this as convenient notation for an ordered
$n$-tuple of weighted points $(x_i,w_i)$.  An $n$-variable \emph{weighted mean} $\gamma$ on a set $X$ assigns to
each weighted $n$-tuple $(\Bw,\mathbf{x})$ some $\gamma(\Bw,\mathbf{x})\in X$ with the extra condition that
$\gamma(\Bw,(x,x,\ldots,x))=x$. We may think of $\gamma(\Bw,\mathbf{x})$ as the assigning of a  ``center of mass."

\begin{remark}  The \emph{weighted geometric mean} $A\#_tB$ (with $\Bw=(1-t,t)$\,) is given by
$$A\#_t B=A^{1/2}(A^{-1/2}BA^{-1/2})^t A^{1/2} .$$
The map $\alpha:[0,\delta(A,B)]\to \PP$ defined by $\alpha(t)=A\#_t B$ is an isometry onto the geodesic arc
between $A$ and $B$.
\end{remark}

Rather obvious variants of the axiomatic properties (P1)-(P10) exist in the weighted mean setting and were introduced
and studied in \cite{LLL}.  There a weighted version of the mean given by Bini, Meini, and Poloni was introduced, which
could also be extended to more general metric spaces.

Closely related to the notion of a omnivariable weighted mean is
that of a barycentric map. Denote by $\mathcal{P}^{<\infty}(X)$ the
set of all finitely supported probability measures on  $X$, measures of the form
$\sum_{k=1}^n w_k\delta_{x_k}$, where $(w_1,\ldots,w_n)$ is a weight
and $\delta_x$ is the unit point mass at $x$.    A \emph{barycentric map}
in its simplest manifestation is a map $\beta:\mathcal{P}^{<\infty}(X)\to X$ such that
$\beta(\delta_x)=x$ for each $x\in X$.  A barycentric $\beta$
gives rise to a corresponding multivariable \emph{weighted mean}
defined by
$$\gamma_\beta(\Bw;x_1,\ldots,x_n)=\beta\left(\sum_{k=1}^n w_k
\delta_{x_k}\right),$$ where $\Bw=(w_1,\ldots,w_n)$ is a weight.  We
note, however, that it is not the case that all weighted means arise
from barycentric maps. If we restrict to uniform weights ($w_i=1/n$
for each $i$), we are essentially in the setting of non-weighted
means.

Soon after the introduction of the ALM mean, a better candidate for the
multivariable matrix geometric mean was put forth to which we now turn.
Let $(M,d)$ be a metric space.  The \emph{least squares mean} $\Lambda(a_1,\ldots,a_n)$
is defined as the solution to the optimization problem of minimizing the sum of distances squared
$$\Lambda(a_1,\ldots,a_n)=\underset{x\in M}{\argmin}\sum_{i=1}^n d^2(x,a_i),$$
and the \emph{weighted least squares
mean} $\Lambda(\Bw;a_1,\ldots,a_n)$ for $\mathbf{w}=(w_1,\ldots, w_n)$ by
$$\Lambda(\Bw;a_1,\ldots,a_n)=\underset{x\in M}{\argmin}\sum_{i=1}^nw_{i} d^2(x,a_i),$$
provided the solution uniquely exists in each of the respective
cases. Note that the least squares mean is equal to the weighted
least squares mean for the uniform weight with all entries $1/n$.
The unique solution of the minimizing problem exists in  Hadamard
spaces \cite[Proposition 1.7]{St03}, since the non-negative function
defined by  $x\mapsto \sum_{i=1}^n w_i d^2(x,a_i)$ is  uniformly
convex in this case. (We recall $F:M\to \R$ is \emph{uniformly
convex} if there exists a strictly increasing
$\varphi: [0,\infty)\to [0,\infty)$ such that for
all $x,y\in M$,
$$F(m(x,y))\leq \frac{1}{2} F(x)+\frac{1}{2}F(y)-\varphi(d(x,y)),$$
where $m(x,y)$ is the unique midpoint between $x$ and $y$.)

E.\ Cartan considered such ``barycenters" in the case of Riemannian manifolds,
where they uniquely exist for the ones of nonpositive curvature and exist locally much
more generally, and M. Fr\'echet considered them in more general metric spaces. Thus the
least squares mean is also called the Cartan mean or Fr\'echet mean.  Such means also
frequently go by the name ``Karcher means," but Karcher's  approach from differential geometry
involved finding the solution of the ``Karcher equation" that was satisfied at this extremum \cite{Kar}.
We return to this in a later section.

M.\ Moakher (2005) \cite{Mo05} first introduced and studied the
 least squares mean for the set of positive definite
matrices $\PP$ equipped with the Riemannian metric as a omnivariable
generalization of the two-variable geometric mean. Independently R.\ Bhatia
and J.\ Holbrook (2006) \cite{BH06,BH06b} introduced and studied the least
squares mean in the weighted setting.  These authors established
its (unique) existence and verified most of the  axiomatic
properties (P1)-(P10) satisfied by the Ando-Li-Mathias geometric
mean: consistency with scalars, joint homogeneity, permutation
invariance, congruence invariance, and self-duality (the last two
being true since congruence transformations and inversion are
isometries). Further, based on computational experimentation, Bhatia
and Holbrook conjectured monotonicity for the least squares mean (property
(P4) in the earlier list), but this was left as an open problem.

\section{The Inductive Mean, Random Variables, and Monotonicity}
One other mean played an important role in what followed, one that
we shall call the \emph{inductive mean}, following the terminology
of K.-T. Sturm (2003) \cite{St03}.  It appeared elsewhere in the
work of M. Sagae and K. Tanabe (1994) \cite{ST} and Ahn, Kim, and
Lim (2007) \cite{AKL}.  It is defined inductively for Hadamard
spaces (or more generally for metric spaces with weighted binary
means $x\#_t y$) for each $k\geq 2$  by $S_2(x,y)=x\# y$ and for
$k\geq 3$, $S_{k}(x_1,\ldots, x_{k})=S_{k-1}(x_1,\ldots,
x_{k-1})\#_{\frac{1}{k}}x_{k}$. (Here $x\#_t y$ is the unique point
$z$ such that $d(x,y)=(1-t)d(x,z)+ td(y,z)$  for $0\leq t\leq 1$.)
Note that this mean at each stage is defined from the previous stage
by taking the appropriate two-variable weighted mean, which is
monotone for the Hadamard space ${\Bbb P}$.  Thus
the inductive mean is monotone (property (P4)).

Let $(M,d)$ be  Hadamard space and $\Bw=(w_1,\ldots, w_n)$ be a weight.
Set  $\N_n=\{1,2,\ldots, n\}$ and assign to $k\in\N_n$ the probability $w_k$.  Let $\Bx=(x_1,\ldots,x_n)\in M^n$.
For each $\omega \in \prod_{j=1}^\infty \N_n$, a countable product, define a sequence $\sigma_{\omega,\Bx}$
in $M$ by $\sigma_{\omega,\Bx}(1)=x_{\omega(1)}$, $\sigma_{\omega,\Bx}(k)=S_k(x_{\omega(1)},\ldots,x_{\omega(k)})$,
where $S_k$ is the inductive mean.   The sequence $\sigma_{\omega,\Bx}$ may be viewed as a ``walk" starting
at $\sigma_{\omega,\Bx}(1)=x_{\omega(1)}$ and moving toward $x_{\omega(k)}$
from $\sigma_{\omega,\Bx}(k-1)$  a distance of
$(1/k)d(\sigma_{\omega,\Bx}(k-1),x_{\omega(k)})$ to obtain obtaining $\sigma_{\omega,\Bx}(k)$. Alternatively we may
give $\prod_{j=1}^\infty \N_n$ the product probability, making it a
probability space, and define a  family of i.i.d.\,random variables $\{X_k\}$ by $X_k(\omega)=x_{\omega(k)}$.  We replace the
traditional sum of the first $k$ random variables by $\sigma_{\omega,\Bx}(k)$ and take for the expected value
$\Lambda(\Bw;x_1,\ldots,x_m)$. From this viewpoint we have the following special case
of Sturm's Law of Large Numbers for Hadamard spaces \cite[Theorem 4.7]{St03}:

\smallskip\noindent
\textbf{Sturm's Theorem}.  \emph{Giving $\prod_{n=1}^\infty
\mathbb{N}_m$ the product probability, the set $$\left\{\omega\in
\prod_{k=1}^\infty \N_m:\lim_{k\to\infty} \sigma_{\omega,\Bx}(k)=
\Lambda(\Bw;x_1,\ldots,x_m)\right\}$$ has measure $1$, i.e.,
$\sigma_{\omega,\Bx}(k)\to \Lambda(\Bw;x_1,\ldots,x_m)$ as
$k\to\infty$ for almost all $\omega$}.

\smallskip We shall briefly
return to the theory of random variables on a probability space
taking values in $\PP$, or more generally in a Hadamard space, at a
later point.

Using the preceding version of Sturm's Theorem, Lawson and Lim (2011)  \cite{LL1} provided a positive solution to the
earlier mentioned conjecture of Bhatia and Holbrook about the monotonicity of the least squares mean.
\begin{theorem}\label{T:mono}
Let $\PP$ be the open cone of positive definite matrices of some fixed dimension, and let $n\geq 3$
\begin{itemize}
\item[(1)] The [weighted] least squares mean $\Lambda$ on $\PP$ is monotone:  $A_i\leq B_i$  for $1\leq i\leq n$ implies $\Lambda(A_1,\ldots, A_n)\leq \Lambda
(B_1,\ldots B_n)$ $[\Lambda(\Bw;A_1,\ldots, A_n)\leq \Lambda
(\Bw;B_1,\ldots B_n)]$.
\item[(2)]  The other nine (weighted) ALM axioms hold for $\Lambda$.
\end{itemize}
\end{theorem}
\textbf{Proof}.  Assume for $\mathbb{A}=(A_1,\ldots,A_m)$ and
$\mathbb{B}=(B_1,\ldots,B_m)$  that
 $A_i\leq B_i$ for $1\leq i\leq
m$. Let $\Bw$ be a weight. By Sturm's theorem applied to the
``walks" for $\mathbb{A}$ and for $\mathbb{B}$, we have
$$\sigma_{\omega,\mathbb{A}}(k)\to \Lambda(\Bw;A_1,\ldots,A_n)~~\mathrm{and}~~
\sigma_{\omega,\mathbb{B}}(k)\to \Lambda(\Bw;B_1,\ldots,B_n)$$
as $k\to\infty$ for almost all  $\omega\in \prod_{n=1}^\infty \N_m$ (since the intersection of two sets of measure 1 has measure 1).
Fixing any such $\omega$,
we obtain part (1) since the partial order relation is closed and by the monotonicity of the inductive mean for each $k$
$$\sigma_{\omega,\mathbb{A}}(k)=S_k(A_{\omega(1)},\ldots,A_{\omega(k)})\leq
S_k(B_{\omega(1)},\ldots,B_{\omega(k)})=\sigma_{\omega,\mathbb{B}}(k).$$

\medskip
Later in \cite{BK12} R.\,Bhatia and
R.\,  Karandikar proved the monotonicity property using in place of Sturm's theorem some elementary counting arguments
and basic inequalities for the metric  $\delta$. Meanwhile J. Holbrook \cite{Hol} identified a specific
$\omega\in \prod_{n=1}^\infty \N_m$ that always yielded a ``walk" converging
 to $\Lambda(A_1,\ldots, A_n)$ for any choice of $A_1,\ldots, A_n\in\PP$. Later Y. Lim and
 M. Palfia \cite{LP14} obtained more general results about deterministic walks yielding the least
 squares mean in general Hadamard spaces.

\smallskip \textbf{Note}: The ALM mean is typically distinct from the least squares mean for $n\geq 3$.  Thus
the ALM axioms do not characterize a mean.  The latter fact had already been noted by
 Bini, Meini and Poloni (2010), who observed their
 variant of the ALM mean was different from it \cite{BMP}. In \cite{LP}, Lim and
 Palfia  showed that the Karcher mean is uniquely determined by
  congruence invariancy (P6), self-duality (P8), and the following Yamazaki inequality \cite{Y}:
  $$\sum_{k=1}^nw_k\log A_k\leq 0\Longrightarrow \Lambda({\bf w};
  A_1,\dots,A_n)\leq I.$$

\section{The Hilbert Space Setting}
Let $\mathcal{B}(E)$ be the $C^*$-algebra of bounded linear
operators equipped with the operator norm on an infinite-dimensional
Hilbert space $E$.  Let $\mathbb{H}(E)$ denote the closed subspace
of hermitian operators, and let $\PP_E$ be the cone of invertible
positive hermitian operators, an open cone in $\mathbb{H}(E)$. Again
$\exp$ and $\log$ are analytic inverses between $\mathbb{H}(E)$ and
$\PP_E$.  We equip $\PP_E$ with the Thompson metric defined by
$d_T(A,B)=\Vert \log(A^{-1/2}BA^{-1/2})\Vert$, where
$\Vert\cdot\Vert$ is the operator norm.  The metric $d_T$ retains
the four properties of $\delta$ in Proposition \ref{P:metprops} in
the finite-dimensional setting, except that $A\#B$ is no longer the unique
midpoint between $A$ and $B$ \cite{Nu}.  Additionally $\PP_E$ is no longer a Hadamard
space, basically because $\mathbb{H}(E)$ is no longer a Hilbert
space.  The definition and basic properties of the geometric mean
$A\# B$, in particular those of Proposition \ref{P:basicprop},
remain valid, and it still gives a distinguished metric midpoint (no
longer unique \cite{Nu}) between $A$ and $B$. The
smallest closed set containing $A$ and $B$ and closed under taking
$\#$-midpoints yields a distinguished metric geodesic connecting $A$
and $B$ consisting of all $A\#_t B$, $0\leq t\leq 1$.

\section{The Karcher Equation}
The uniform convexity of the Riemannian metric $\delta$  on $\PP$
yields that the least squares mean is the unique critical point for
the function $X\mapsto \sum_{k=1}^n
w_k\delta^2(X,A_k)$. The least squares mean is
thus characterized by the vanishing of the gradient, which is
equivalent to its being a solution of the following \emph{Karcher
equation}:
\begin{equation}\label{E:Karch}
\sum_{k=1}^{n}w_{k}\log (X^{-1/2}A_{k}X^{-1/2})=0.
\end{equation}

The Karcher equation (\ref{E:Karch})
can also be used to \emph{define}  a mean on the cone $\PP$ of positive  invertible  bounded operators on
an infinite-dimensional Hilbert space (where one no longer has a Hadamard space), called the \emph{Karcher mean}.
As we just previously noted, restricted to the matrix setting it yields the least squares mean. Thus it is reasonable
to continue to denote it by $\Lambda(\Bw,A_1,\ldots,A_n)$.

Power means for positive definite matrices were introduced
by Lim and Palfia (2012) \cite{LP}.

\smallskip\noindent
\textbf{Theorem}.  \emph{Let $A_{1},\dots,A_{n}\in\PP$ and let
$\Bw=(w_{1},\dots,w_{n})$ be a weight.  Then for each  $ t\in
(0,1],$  the following equation has a unique positive definite
solution} $X=P_t(\Bw;A_1,\ldots,A_n)$,  \emph{called the $t$-weighted power mean}:
\begin{eqnarray*}
X=\sum_{k=1}^{n}w_{k}(X\#_{t}A_{i}).
\end{eqnarray*}

\medskip When restricted to the positive reals,
the power mean reduces to the usual power mean
$$P_t(\Bw;a_1,\ldots,a_n)=\left(w_{1}a_{1}^{t}+\cdots+w_{n}a_{n}^{t}\right)^{\frac{1}{t}}.$$

In 2014 Lawson and Lim showed  that the preceding notion of power mean extended to the setting of
bounded operators on a Hilbert space  \cite{LL14}
and  established that the power means are decreasing,
$s<t$ implies $P_s(\cdot\,;\,\cdot)\leq P_t(\cdot\,;\,\cdot)$.
Using power means Lawson and Lim were able to establish the existence and uniqueness of the Karcher mean
in the $C^*$-algebra of bounded operators on a Hilbert space.

\begin{theorem}\label{T:karch}
 In the strong operator topology
 $$\Lambda(\cdot\,;\,\cdot)=\lim_{t\to 0^+} P_t(\cdot\,;\,\cdot)=\inf_{t>0}P_t(\cdot\,;\,\cdot),$$
where $\Lambda$ is the Karcher mean, the unique solution of the Karcher equation
$$X=\Lambda({\bf w};A_1,\ldots,A_n)\Leftrightarrow \sum_{k=1}^n w_k\log(X^{-1/2}A_kX^{-1/2})=0.$$
\end{theorem}

Via this machinery many of the axiomatic properties of the least squares mean in the finite-dimensional
setting were extended to the corresponding Karcher mean in the infinite-dimensional setting.

Recent work by Lim-Palfia \cite{LP18} and independently by Lawson
\cite{Law} shows that the preceding constructions and results remain
valid  for the open cone of positive invertible elements in any
unital $C^*$-algebra.

\section{Barycenters}
A \emph{Borel probability measure} on a metric space $(X,d)$ is a
countably additive non-negative measure $\mu$ on the Borel algebra
$\mathcal{B}(X)$, the smallest $\Sigma$-algebra containing the open
sets, such that $\mu(X)=1$.   We denote the set of all probability
measures on $(X,\mathcal{B}(X))$ by $\Pro(X)$. Let ${\mathcal
P}_{0}(X)$ be  the set of all uniform finitely supported probability
measures, i.e., all $\mu\in {\mathcal P}(X)$ of the form
$\mu=\frac{1}{n}\sum_{k=1}^{n}\delta_{x_{k}}$ for
some $n\in{\mathbb N},$ where $\delta_x$ is the point measure of
mass $1$ at $x$.

A measure $\mu\in\Pro(X)$ is said to be \emph{integrable}  if for some (and hence all) $x$
$$\int_{X}d(x,y)d\mu(y)<\infty.$$
The set of integrable measures is denoted by $\Pro^1(X)$.


The \emph{Wasserstein distance}  (alternatively Kantorovich-Rubinstein distance) $d^W$
on ${\mathcal P}^{1}(X)$ is a standard metric for probability measures.
 It is known that $d^W$  is a complete metric on ${\mathcal P}^{1}(X)$ whenever $X$ is a complete metric space
and that ${\mathcal P}_{0}(X)$ is $d^W$-dense in ${\mathcal P}^{1}(X)$.
\begin{definition}(\cite{St03})
A map $\beta:(\Pro^1(X),d^W)\to (X,d)$ is a \emph{contractive
barycentric map} if (i) $\beta(\delta_x)=x$ for all $x\in X$, and
(ii) $d(\beta(\mu),\beta(\nu))\leq d^W(\mu,\nu)$ for all $\mu,\nu\in
\Pro^1(X)$.
\end{definition}

Suppose $G=\{G_n:X^n\to X: n\geq 2\}$ is a omnivariable mean on $(X,d)$, a metric space.  The mean $G$
is said to be \emph{iterative} if for all $n,k\geq 2$ and $\Bx=(x_1,\ldots,x_n)\in X^n$,
$$G_n(\Bx)=G_{nk}(\Bx^k),$$
where $\Bx^k=(x_1,\ldots,x_n, x_1,\ldots,x_n,\ldots, x_1,\ldots,x_n)\in X^{nk}$.  The next proposition is a key result
from \cite{LL16}.
\begin{proposition}\label{P:bary}
Suppose $G$ is a omnivariable mean on a complete metric space $(X,d)$.   If $G$ is symmetric, iterative,
and satisfies for all $n$, all $\Bx=(x_1,\ldots,x_n)$, $\By=(y_1,\ldots,y_n)$,
\begin{equation}\label{E:Wasser}
 d(G(\Bx),G(\By))\leq \frac{1}{n}\sum_{k=1}^n d(x_k,y_k), \end{equation}
then there exists a unique contractive barycentric map
$\beta:\Pro^1(X)\to X$ satisfying $G(x_1,\ldots,
x_n)=\beta\left(\sum_{k=1}^n
(1/n)\delta_{x_k}\right)$.
\end{proposition}
The Wasserstein distance between $\sum_{k=1}^n \frac{1}{n}\delta_{x_k}$
and $\sum_{k=1}^n \frac{1}{n}\delta_{y_k}$ is the minimum
over all permutations of $\{y_k\}_{k=1}^n$ of the sum on the right in equation (\ref{E:Wasser}), so the inequality
implies $\beta$ is contractive on $\Pro_0(X)$, a dense subset, and hence extends uniquely to a contractive
map on $\Pro^1(X)$.

The hypotheses of the preceding proposition are satisfied by the Karcher mean $\Lambda$, and thus the
mean ``extends" to a contractive barycentric map $\beta_\Lambda:\Pro^1(\PP)\to \PP$, called the
\emph{Karcher barycentric map}.  It is characterized by
$$ X=\beta_\Lambda(\mu) \Leftrightarrow \int_\PP \log(X^{-1/2} A X^{-1/2}) d\mu(A)=0.$$
The existence and basic theory and properties of the Karcher
barycentric map can be found in \cite{LL16}. Many of the basic
properties of the least squares mean in Theorem \ref{T:ALM} have
analogous properties for the Karcher barycentric map.

\section{$\PP$-valued Random Variables}
The Karcher barycentric map makes possible a theory of $\PP$-valued
random variables.  Let $(Z,\mathcal{A},q)$ be a probability space,
where $Z$ is a set, $\mathcal{A}$ is a $\sigma$-algebra, and $q$ is
a probability measure on $\mathcal{A}$.  Let $X:Z\to\PP$ be
measurable. With respect to $X$ there is \emph{push-forward} measure
$X_*(q)$ on $\PP$ defined by $X_*(q) (Q)=q(X^{-1}(Q))$ for any Borel
subset $Q$ of $\PP$. This measure is called the \emph{distribution}
of $X$.  We are particularly interested in the $L^1$ random
variables, those for which $X_*(q)\in
\Pro^1(\PP)$. In this case we can define the \emph{expectation}
$E(X)=\beta(X_*(q))$, where $\beta$ is the Karcher barycentric map.
These ideas can be worked out in Hadamard spaces to a rather
full-blown theory as has been done by Sturm in \cite{St03}, and we
refer the reader to that source for details.  However, the infinite
dimensional $\PP$ is not a Hadamard space, so
requires new approaches; see \cite{LP18}.

\section{Summary}
In the preceding we have attempted to outline the high points of the striking development of the theory of the matrix/operator geometric mean
on the cone of positive matrices/operators in the past twenty years. Over this period of time it has evolved from  a two-variable matrix mean
to a omnivariable matrix mean
(the least squares, Cartan, or Frech\'{e}t mean) to an operator
mean in the setting of unital $C^*$-algebras (the Karcher mean) and finally to a barycentric map on the space of integrable Borel probability measures.
At each stage of the evolution significant new insights and developments were necessary. The theory has drawn heavily from matrix and operator
theory and at the same time from geometric notions.
And along the way we have seen a variety of characterizations
of this mean: the least squares mean, the probabilistic or deterministic characterization as a limit of a ``walk" with the inductive mean,
the solution of the Karcher equation, and the limit/infimum of the power means $\{P_t\}$ for $t\searrow 0$.
Whatever future developments may hold, it is clear that a substantial theory has already emerged.

\end{document}